\documentclass{article}
\usepackage{amssymb}
\usepackage{amsfonts}
\usepackage{graphicx}
\usepackage{amsmath}
\usepackage{booktabs}
\usepackage{multirow}

\newtheorem{algorithm}{Algorithm}

\newtheorem{proposition}{Proposition}

\newenvironment{proof}[1][Proof]{\textbf{#1.} }{\ \rule{0.5em}{0.5em}}

\begin{document}

\date{}
\title{Gaussian rule for integrals involving Bessel functions}
\author{Eleonora Denich\thanks{%
Dipartimento di Matematica e Geoscienze, Universit\`{a} di Trieste,
eleonora.denich@phd.units.it} \and Paolo Novati \thanks{%
Dipartimento di Matematica e Geoscienze, Universit\`{a} di Trieste,
novati@units.it}}
\maketitle

\begin{abstract}
In this work we develop the Gaussian quadrature rule for weight functions involving fractional powers, exponentials and Bessel functions of the first kind. Besides the computation based on the use of the standard and the modified Chebyshev algorithm, here we present a very stable algorithm based on the preconditioning of the moment matrix. Numerical experiments are provided and a geophysical application is considered.
\end{abstract}



\section{Introduction}

This work deals with the construction of Gaussian quadrature rule for the computation of integrals of the type
\begin{equation} \label{integrale}
I_{\nu,\alpha,c}(f) = \int_0^{\infty} f(x) x^{\alpha} e^{-cx} J_{\nu} (x) dx,
\end{equation}
where $J_{\nu}$ is the Bessel function of the first kind of order $\nu \geq 0$, $\alpha>-1$, $c>0$ and $f$ is a smooth function.
Since for the Bessel functions it holds $|J_{\nu}(x)| \leq 1$, for $\nu \geq 0$, $x \in \mathbb{R}$ (see \cite[p.362]{Abramowitz}), we consider weight functions of the type
\begin{equation} \label{peso}
w_{\nu, \alpha,c}(x)=x^{\alpha} e^{-cx} [J_{\nu}(x)+1] \quad {\rm on} \; [0, + \infty).
\end{equation}
Then, we rewrite (\ref{integrale}) as
\begin{equation*} 
I_{\nu,\alpha,c}^J(f)-I_{\alpha,c}^L(f),
\end{equation*}
where
\begin{equation} \label{integrale2}
I_{\nu,\alpha,c}^J(f)=\int_0^{\infty} f(x) x^{\alpha} e^{-cx} [J_{\nu} (x)+1] dx,
\end{equation}
and
\begin{equation} \label{integrale3.1}
I_{\alpha,c}^L(f)=\int_0^{\infty} f(x) x^{\alpha} e^{-cx} dx.
\end{equation}
We notice that the integral (\ref{integrale3.1}) can be accurately computed using a slight modification of the Gauss-Laguerre quadrature rule.
In this setting, our aim is to construct a Gaussian rule with respect to the function $w_{\nu, \alpha,c}$. Since we do not know the explicit expression of the corresponding monic orthogonal polynomials, that we denote by $\pi_k$, $k \geq 0$, we need to employ a numerical scheme to derive the coefficients of the three-term recursion
\begin{equation*}
\pi_{k+1}(x) = (x- \alpha_k) \pi_k(x) - \beta_k \pi_{k-1}(x), \quad k \geq 0,
\end{equation*}
\begin{equation*}
\pi_{-1}(x) = 0, \quad \pi_0(x)=1,
\end{equation*}
with $\beta_k > 0$.
This can be done by computing the associated moments 
\begin{equation} \label{PM}
\mu_k^{\nu, \alpha,c} = \int_0^{\infty} x^k w_{\nu,\alpha,c} (x) dx, \quad k \geq 0,
\end{equation}
and then using the Chebyshev algorithm (see \cite[sect.2.3]{Gautschi1982}).
These coefficients define the tridiagonal symmetric Jacobi matrix, whose eigenvalue decomposition provides abscissas and weights of the quadrature rule.
This final step is efficiently implemented by the famous Golub and Welsh algorithm \cite{Golub}. Some alternatives to this algorithm have been later developed and we refer to \cite{Laurie} for a general discussion and a rich bibliography.
Nevertheless, it is well known (see e.g. \cite{Gautschi}) that the computation of the recurrence coefficients can be inaccurate for growing $k$ because  the problem is severely ill conditioned when starting from the power moments (\ref{PM}).
The ill conditioning can be partially overcome by using the modified moments, having at disposal a family of polynomials orthogonal with respect to a weight function similar to the one of the problem. This approach may be efficient in general but not always when working with unbounded intervals of integration (see \cite{Gautschi726} and \cite{Gautschi}). The idea of using modified moments was introduced by Sack and Donovan in \cite{SD}, who developed an algorithm similar to the so called modified Chebyshev algorithm, advanced by Gautschi in \cite[sect.2.4]{Gautschi1982}. The same algorithm was independently obtained by Wheeler in \cite{Wheeler}.

In this work we present an alternative approach that is based on the preconditioning of the moment matrix. In particular, since the three-term recurrence coefficients can be written in terms of ratios of determinants of the moment matrix or slight modifications of them (see \cite[sect.2.7]{Davis}), we exploit the Cramer rule to show that the coefficients can be computed by solving a linear system with the moment matrix. Since the weight function (\ref{peso}) can be interpreted as a perturbation of the weight function of the generalized Laguerre polynomials, we use the moment matrix of these polynomials as preconditioner.  
The numerical experiments show that this technique is always (independently of the parameters $\nu, \alpha,c$) much more stable than the modified Chebyshev algorithm.

As an application we use the developed Gaussian quadrature to evaluate integrals of the type (\ref{integrale}) arising in geophysical electromagnetic (EM) survey. In particular, we consider the electromagnetic fields  over a layered earth due to magnetic dipoles above the surface (see \cite{WH}). In this framework, $f$ is a smooth function, $\nu=0,1$, $\alpha=0$ and $0<c<1$.

The paper is organized as follows. In Section \ref{section2} we derive a recursive relation for the practical evaluation of the power moments. In Section \ref{section3} we show the necessary details for the construction of the Gaussian rule for (\ref{integrale}) by using the Chebyshev algorithm. In Section \ref{section4} we employ the modified Chebyshev algorithm working with the modified moments generated by the generalized Laguerre polynomials. In Section \ref{section5} we present the alternative approach based on the preconditioning of the moment matrix, using again the generalized Laguerre polynomials. Finally, in Section \ref{section6} we apply the method for the computation of EM fields. 

\section{Computation of the moments} \label{section2}

In order to compute the moments
\begin{equation} \label{mom}
\mu_k=\mu_k^{\nu, \alpha ,c} = \int_0^{\infty} x^{k+\alpha} e^{-cx} [J_{\nu}(x)+1] d x, \quad k \geq 0,
\end{equation}
we first derive a recursive relation for the so called core moments, defined as
\begin{equation} \label{defCORE}
\mu_{k,0}=\mu_{k,0}^{\nu,\alpha ,c} = \int_0^{\infty} x^{k+\alpha} e^{-cx} J_{\nu}(x) d x, \quad k \geq 0.
\end{equation}

\begin{proposition}
For $k \geq 0$ it holds
\begin{equation} \label{core}
\mu_{k,0} = \frac{1}{(\sqrt{c^2+1})^{k+\alpha+1}} \Gamma (k+\alpha+ \nu +1) P_{k+\alpha}^{-\nu} \left( \frac{c}{\sqrt{c^2+1}} \right),
\end{equation}
where $\Gamma$ is the Gamma function and $ P_{k+\alpha}^{- \nu}$ is the associated Legendre function (see e.g. \cite[ch.8]{Abramowitz}) of order $-\nu$ and degree $k+\alpha$.
\end{proposition}
\begin{proof}
We start from the general relation \cite[p.713]{GR}
\begin{equation*}
\int_0^{\infty} e^{-t \cos \theta} J_{\nu} (t \sin \theta) t^{k+\alpha} d t = \Gamma(k+ \alpha+ \nu +1) P_{k+\alpha}^{- \nu} (\cos \theta),
\end{equation*}
which holds for each $k \geq 0$ whenever $\alpha > -1$, $\nu \geq 0$.
By the change of variable $s=t \sin \theta$, we have that
\begin{equation*}
\int_0^{\infty} e^{-s \frac{\cos \theta}{\sin \theta}} J_{\nu} (s) s^{k+\alpha} d s = \sin^{k+\alpha+1}(\theta) \Gamma(k+ \alpha+ \nu +1)P_{k+\alpha}^{- \nu} (\cos \theta).
\end{equation*}
Setting $\theta=\arctan \left( \frac{1}{c} \right)$, $0 < \theta <\frac{\pi}{2}$, so that $c=\frac{\cos \theta}{\sin \theta}$, and using the relations
\begin{equation*}
\sin (\arctan x)= \frac{x}{\sqrt{1+x^2}}, \quad \cos (\arctan x)= \frac{1}{\sqrt{1+x^2}},
\end{equation*}
we obtain the result.
\end{proof}

\begin{proposition}
The following three-term recursion holds
\begin{equation} \label{ric_core}
\mu_{k+1,0} = \frac{1}{c^2+1}\lbrace c \left[ 2(k+\alpha)+1 \right] \mu_{k,0} -[(k+\alpha)^2- \nu^2]\mu_{k-1,0} \rbrace, \quad k \geq 1,
\end{equation}
with
\begin{align}
\mu_{0,0}&= \frac{\Gamma(\alpha+\nu+1)}{(\sqrt{c^2+1})^{\alpha+1} \Gamma(\nu+1)} \left( \frac{c+\sqrt{c^2+1}}{\sqrt{c^2+1}-c} \right)^{-\frac{\nu}{2}} \times \notag \\ 
& \quad {}_2F_1 \left(-\alpha, \alpha+1; 1+\nu; \frac{\sqrt{c^2+1}-c}{2 \sqrt{c^2+1}} \right), \label{core0}\\
\mu_{1,0}&= \frac{\Gamma(\alpha+\nu+2)}{(\sqrt{c^2+1})^{\alpha+2} \Gamma(\nu+1)} \left( \frac{c+\sqrt{c^2+1}}{\sqrt{c^2+1}-c} \right)^{-\frac{\nu}{2}} \times \notag \\
& \quad {}_2F_1 \left(-\alpha-1, \alpha+2; 1+\nu; \frac{\sqrt{c^2+1}-c}{2 \sqrt{c^2+1}} \right), \label{core1}
\end{align}
where ${}_2F_1$ is the hypergeometric function.
\end{proposition}
\begin{proof}
From equation (\ref{core}) and using the following three-term recursive relation for the associated Legendre functions (\cite[p.334]{Abramowitz}) 
\begin{equation*}
(k+\alpha+\nu+1) P_{k+\alpha+1}^{-\nu} (z) = (2k+2\alpha+1) z P_{k+\alpha}^{-\nu} (z) - (k+\alpha-\nu) P_{k+\alpha-1}^{-\nu}(z),
\end{equation*}
we can write
\begin{align}
\mu_{k+1,0} &= \frac{1}{(\sqrt{c^2+1})^{k+\alpha+2}} \Gamma (k+ \alpha + \nu +2) P_{k+\alpha+1}^{-\nu} \left( \frac{c}{\sqrt{c^2+1}} \right) \notag \\
&= \frac{\Gamma(k+\alpha+\nu+2)}{(\sqrt{c^2+1})^{k+\alpha+2}} \Bigg[ \frac{2(k+\alpha)+1}{k+\alpha+\nu+1} \frac{c}{\sqrt{c^2+1}} P_{k+\alpha}^{-\nu} \left( \frac{c}{\sqrt{c^2+1}} \right) \notag \\
&- \frac{k+\alpha- \nu}{k+\alpha+ \nu+1} P_{k+\alpha-1}^{- \nu} \left( \frac{c}{\sqrt{c^2+1}} \right) \Bigg]. \label{ricLeg}
\end{align}
Rearranging (\ref{ricLeg}) and using again (\ref{core}) for $\mu_{k,0}$ and $\mu_{k-1,0}$, we obtain the relation (\ref{ric_core}). Equations (\ref{core0}) and (\ref{core1}) follow directly from (\ref{core}) with $k=0$ and $k=1$, respectively, and from the relation \cite[p.999]{GR}
\begin{equation*}
P_{k+\alpha}^{-\nu} (x)=\frac{1}{\Gamma(\nu+1)} \left(\frac{1+x}{1-x} \right)^{-\frac{\nu}{2}} {}_2F_1 \left(-k-\alpha,k+\alpha+1;1+\nu;\frac{1-x}{2} \right),  
\end{equation*}
for $x \in (0,1)$.
\end{proof}

Finally, we can derive a recursive relation for the moments. 
\begin{proposition}
For $k \geq 1$ it holds
\begin{align} 
\mu_{k+1} &= \frac{1}{c^2+1} \Bigg\lbrace c \left[ 2(k+\alpha)+1 \right] \mu_k -\left[(k+\alpha)^2-\nu^2\right] \mu_{k-1} \notag \\ 
&+  \frac{\Gamma(k+\alpha)[(k+\alpha)^2+(k+\alpha)-c^2\nu^2]}{c^{k+\alpha+2}} \Bigg\rbrace, \label{ricorsione_momenti}
\end{align}
with
\begin{equation*}
\mu_0 = \mu_{0,0}+\frac{\Gamma(\alpha+1)}{c^{\alpha+1}}, \quad \mu_1 = \mu_{1,0}+\frac{\Gamma(\alpha+2)}{c^{\alpha+2}}.
\end{equation*}
\end{proposition}
\begin{proof}
By definition (\ref{mom}), the moments $\mu_k$ are given by
\begin{align*}
\mu_k&=\int_0^{\infty} x^{k+\alpha} e^{-c x} J_{\nu}(x) dx + \int_0^{\infty} x^{k+\alpha} e^{-c x} dx \\
&=\mu_{k,0} + \frac{\Gamma(k+\alpha+1)}{c^{k+\alpha+1}},
\end{align*}
where we have used \cite[sect.3.381, n.4]{GR}. Therefore, from relation (\ref{ric_core}) for $k \geq 1$, we can write
\begin{align*}
\mu_{k+1} &= \mu_{k+1,0}+\frac{\Gamma(k+\alpha+2)}{c^{k+\alpha+2}} \\
&= \frac{c[2(k+\alpha)+1]}{c^2+1} \left( \mu_k - \frac{\Gamma(k+\alpha+1)}{c^{k+\alpha+1}}\right) - \frac{(k+\alpha)^2-\nu^2}{c^2+1} \left(\mu_{k-1} -\frac{\Gamma(k+\alpha)}{c^k} \right) \\ 
&+ \frac{\Gamma(k+\alpha+1)}{c^{k+\alpha+1}}.
\end{align*}
After some simple manipulations, we obtain the result.
\end{proof}

\section{Computing the three-term recursion} \label{section3}

One of the most used method for the computation of the coefficients $\alpha_k$ and $\beta_k$ of the recurrence relation 
\begin{equation} \label{doppio_asterisco}
\pi_{k+1}(x) = (x- \alpha_k) \pi_k(x) - \beta_k \pi_{k-1}(x), \quad k \geq 0
\end{equation}
\begin{equation*}
\pi_{-1}(x) = 0, \quad \pi_0(x)=1,
\end{equation*}
with $\beta_k >0$, is the Chebyshev algorithm (see \cite[sect.2.3]{Gautschi1982} and \cite{Gautschi}). 

Given the first $2n$ moments $\mu_0,\ldots, \mu_{2n-1}$, the algorithm uniquely determines the first $n$ recurrence coefficients $\alpha_k$ and $\beta_k$, $k=0,\ldots,n-1$, by using the mixed moments
\begin{equation*}
\sigma_{kl} = \int_0^{\infty} \pi_k(x) x^l w_{\nu,\alpha,c} (x) dx, \quad k,l \geq -1.
\end{equation*} 

The Chebyshev algorithm is summarized in Algorithm \ref{ChebyAlg}.

\begin{algorithm} \label{ChebyAlg}
Initialization
\begin{align*}
\alpha_0&=\frac{\mu_1}{\mu_0}, \; \beta_0=\mu_0, \\
\sigma_{-1,l}&=0, \quad l=1,2,\ldots,2n-2, \\
\sigma_{0,l}&= \mu_l, \quad 0,1,\ldots,2n-1,
\end{align*}
for $k=1,2,\ldots,n-1$
\medskip

for $\quad l=k,k+1,\ldots,2n-k-1$
\begin{align*}
\sigma_{k,l}&=\sigma_{k-1,l+1}-\alpha_{k-1} \sigma_{k-1,l}-\beta_{k-1} \sigma_{k-2,l},  \\
\alpha_k &= \frac{\sigma_{k,k+1}}{\sigma_{k,k}}-\frac{\sigma_{k-1,k}}{\sigma_{k-1,l-1}}, \quad \beta_k=\frac{\sigma_{k,k}}{\sigma_{k-1,k-1}}.
\end{align*}
\end{algorithm}

The corresponding Jacobi matrix 
\[ J=
\begin{bmatrix}
\alpha_0 & \sqrt{\beta_1} & & & 0\\
\sqrt{\beta_1} & \alpha_1 & \sqrt{\beta_2} \\
 & \sqrt{\beta_2} & \alpha_2 & \ddots \\
 & & \ddots & \ddots & \sqrt{\beta_{n-1}}\\
0 & & & \sqrt{\beta_{n-1}} & \alpha_{n-1}
\end{bmatrix}  \in \mathbb{R}^{n \times n},
\]
contains the coefficients of the three term recurrence relation for the orthonormal polynomials, that is,
\begin{equation*}
\sqrt{\beta_{k+1}}\tilde{\pi}_{k+1}(x) = (x- \alpha_k) \tilde{\pi}_k(x) - \sqrt{\beta_k} \tilde{\pi}_{k-1}(x), \quad k \geq 0,
\end{equation*}
\begin{equation*}
\tilde{\pi}_{-1}(x) = 0, \quad \tilde{\pi}_0(x)=\frac{1}{\sqrt{\beta_0}}.
\end{equation*}
It is well known that the eigendecomposition of the matrix $J$ provides the nodes $x_i$ and the weights $w_i$, $i=1,\ldots,n$ of the $n$-point Gaussian rule (see e.g. \cite[sect.2.7]{Davis}, and the reference therein).

Finally, for the computation of (\ref{integrale}), we use the approximation
\begin{equation*}
I_{\nu,\alpha,c}^J(f) \approx I_n^J(f) =\sum_{i=1}^n w_i f(x_i),
\end{equation*}
for the integral (\ref{integrale2}). Then, denoting by $t_i^L$, $w_i^L$ respectively the nodes and the weights of the Gauss-Laguerre rule with respect to the weight function $w_{\alpha}(t)=t^{\alpha} e^{-t}$, $\alpha > -1$, the integral (\ref{integrale3.1}) is approximated by
\begin{align*}
I_{\alpha,c}^L(f) &= \frac{1}{c^{\alpha+1}} \int_0^{\infty} f \left( \frac{t}{c} \right) t^{\alpha} e^{-t}  dt \\
& \approx I_n^L(f) = \frac{1}{c^{\alpha+1}} \sum_{i=1}^n w_i^L f \left( \frac{t_i^L}{c} \right).
\end{align*}
Finally, we thus have
\begin{equation} \label{asterisco}
I_{\nu,\alpha,c}(f) \approx I_n^J(f)-I_n^L(f).
\end{equation}

\begin{figure}
\begin{center}
\includegraphics[scale=0.35]{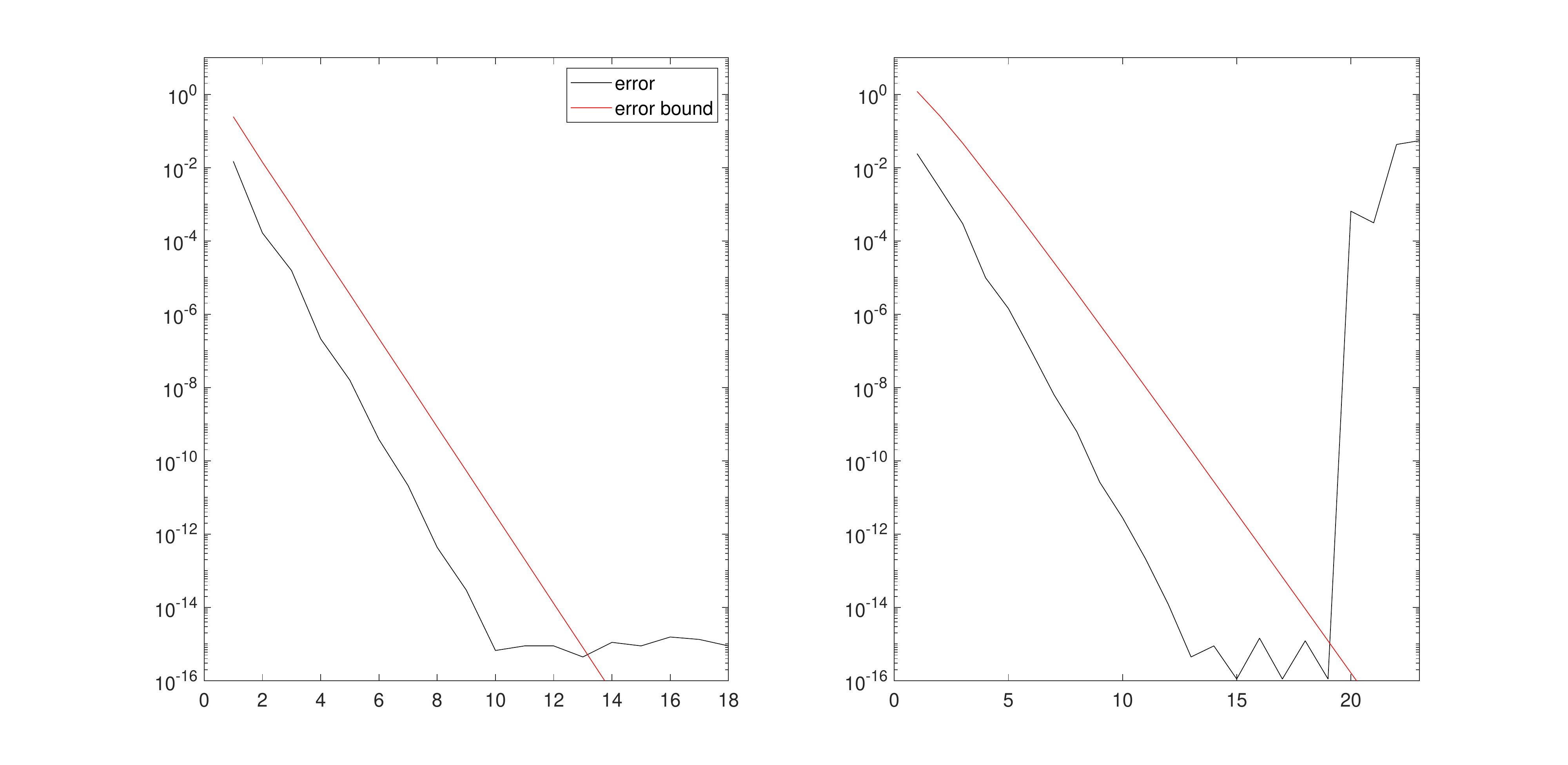}
\end{center}
\caption{Error behavior of (\ref{asterisco}) and error bound (\ref{upper_bound}) with respect to $n$ for $\nu =1$, $\alpha=-0.5$, $c=1$ on the left and for $\nu =0.5$, $\alpha=0.5$, $c=0.2$ on the right. In both cases $f(x)=e^{-0.5x}$.}
\label{esempiCheby}
\end{figure}

Below we present the results of the numerical experiments carried out in Matlab. In particular, the Matlab routine that implements the Chebyshev algorithm is taken from \cite{Gautschi_bis}. Since for integrals involving Bessel functions, exponentials and powers the exact solution is known, in our simulations we choose $f(x) = e^{-0.5 x}$. In Figure \ref{esempiCheby} we consider two examples, for different values of the parameters $\nu,\alpha$ and $c$, and plot the absolute error between the approximation obtained with the developed Gaussian rule and the exact solution (see \cite[sect.6.624, n.6]{GR} and \cite[sect.8.704]{GR}) given by
\begin{align*}
I_{\nu,\alpha,d}(f) &= \frac{\Gamma(d+\nu+1)}{(\sqrt{c^2+1})^{d+1}\Gamma(\nu+1)} \left(\frac{\sqrt{c^2+1}+c}{\sqrt{c^2+1}-c} \right)^{-\frac{\nu}{2}} \times \\
& \quad {}_2F_1 \left( -d,d+1;1+\nu;\frac{1}{2}-\frac{c}{2\sqrt{c^2+1}} \right),
\end{align*}
where $d=c+0.5$.

Moreover, since for the truncation error it holds (see \cite[sect.4.4]{Davis})
\begin{align*} 
E_n(f) &= I_{\nu,\alpha,c}(f)-(I_n^J(f)-I_n^L(f)) \\
&= (I_{\nu,\alpha,c}^J(f)-I_n^J(f))-(I_{\alpha,c}^L(f))-I_n^L(f)) \\
&=\frac{f^{(2n)}(\eta^J)}{(2n)! (k_n^J)^2}-\frac{f^{(2n)}(\eta^L)}{(2n)! (k_n^L)^2}, \quad \eta^J,\eta^L \in (0,\infty),
\end{align*}
where $k_n^J$ and $k_n^L$ are the leading coefficients of the corresponding orthonormal polynomials of degree $n$, in Figure \ref{esempiCheby} we also provide the plot of the upper bound of $E_n(f)$ given by
\begin{equation} \label{upper_bound}
| E_n(f)| \leq \frac{\| f^{(2n)}\|_{\infty}}{(2n)!} \left( \frac{1}{(k_n^J)^2} +\frac{1}{(k_n^L)^2} \right).
\end{equation}
The coefficients $k_n^J$ are numerically evaluated by using the relation (see \cite[sect.2.7]{Davis})
\begin{equation*}
k_n^J = \frac{1}{\prod_{j=0}^n \sqrt{\beta_j}}, 
\end{equation*}
while for $k_n^L$ we employ the known explicit formulation
\begin{equation*}
\quad k_n^L = \frac{1}{\sqrt{n! \Gamma(n+\alpha+1)}}.
\end{equation*}

\section{The modified Chebyshev algorithm} \label{section4}

The picture on the right of Figure \ref{esempiCheby} shows the stability problem when working with the power moments $\mu_k$. Indeed the Chebyshev algorithm typically starts to produce negative values of $\beta_k$ for $k$ around $20$ or even before. This behavior is rather common and has been observed by many authors in the past (\cite{Gautschi}, \cite{Golub}, \cite{Wheeler}).
As already mentioned, the problem is that the coefficients $\alpha_k$ and $\beta_k$ are extremely sensitive to small changes in the moments. In fact, the nonlinear map
\begin{align*}
K_n \colon \mathbb{R}^{2n} &\rightarrow \mathbb{R}^{2n} \\
\mu &\mapsto \rho
\end{align*}
which maps the moment vector $\mu=[\mu_0, \mu_1, \ldots, \mu_{2n-1}]^T$ to the vector $\rho = [\alpha_0,\ldots, \alpha_{n-1}, \beta_0, \ldots, \beta_{n-1}]^T$ of recursion coefficients becomes extremely ill conditioned as $n$ increases (see \cite{Gautschi} for the complete analysis). 

In order to overcome this difficulty, the modified Chebyshev algorithm (see \cite[sect.2.4]{Gautschi1982}) can be employed.
It is based on the use of the modified moments
\begin{equation} \label{modified}
m_k=m_k^{\nu, \alpha,c} = \int_0^{\infty} p_k(x) w_{\nu,\alpha,c} (x)  dx, \quad k \geq 0,
\end{equation}
and on the mixed moments
\begin{equation*}
\tilde{\sigma}_{kl} = \int_0^{\infty} \pi_k(x) p_l(x) w_{\nu,\alpha,c} (x)  dx, \quad k,l \geq -1,
\end{equation*}
where $\lbrace p_k \rbrace_{k \geq 0}$ is a given system of orthogonal polynomials, chosen to be close to the desired polynomials $\lbrace \pi_k \rbrace_{k \geq 0}$, which satisfies the three-term recurrence relation 
\begin{equation*}
p_{k+1}(x) = (x- a_k) p_k(x) - b_k p_{k-1}(x), \quad k \geq 0
\end{equation*}
\begin{equation*}
p_{-1}(x) = 0, \quad p_0(x)=1,
\end{equation*}
with coefficients $a_k \in \mathbb{R}$, $b_k \geq 0$ that are known.

In our case, since the weight function can be interpreted as a perturbation of the weight function relative to the generalized Laguerre polynomials $\lbrace L_k^{\alpha} \rbrace_{k \geq 0}$, we choose as $\lbrace p_k \rbrace_{k \geq 0}$ the system $\lbrace L^{\alpha,c}_k \rbrace_{k \geq 0}$ of the monic polynomials 
\begin{equation} \label{Lac}
L_k^{\alpha,c} (x) = \frac{1}{c^k} \tilde{L}_k^{\alpha} (cx),
\end{equation}
where $\tilde{L}_k^{\alpha} (t)= (-1)^k k! L_k^{\alpha} (t)$ is the monic generalized Laguerre polynomial of degree $k$.
This system satisfies the relation
\begin{equation*}
\tilde{L}_{k+1}^{\alpha} (t)=\left( t-A_k\right) \tilde{L}_{k}^{\alpha}(t)-B_k \tilde{L}_{k-1}^{\alpha} (t), 
\end{equation*}
with
\begin{equation} \label{A_k B_k}
A_k = 2k+\alpha+1, \quad B_k = k(k+\alpha) \quad k \geq 1.
\end{equation}

\begin{proposition}
The monic polynomials $\lbrace L^{\alpha,c}_k \rbrace_{k \geq 0}$ defined by (\ref{Lac}) are orthogonal with respect to the weight function $x^{\alpha} e^{-cx}$ and satisfy the three-term recurrence relation
\begin{equation} \label{Lac_ric}
L_{k+1}^{\alpha,c} (x)=\left( x-\frac{A_k}{c} \right) L_{k}^{\alpha,c}(x)-\frac{B_k}{c^2} L_{k-1}^{\alpha,c} (x), 
\end{equation}
in which $A_k$ and $B_k$ are defined in (\ref{A_k B_k}).
\end{proposition}
\begin{proof}
The orthogonality follows from the change of variable $cx=t$, that leads to
\begin{align*}
\int_0^{\infty} L_k^{\alpha,c} (x) L_l^{\alpha,c} (x) x^{\alpha} e^{-cx} dx &= \frac{1}{c^{k+l}} (-1)^{k+l} k! l! \int_0^{\infty} L_k^{\alpha}(cx) L_l^{\alpha}(cx) x^{\alpha} e^{-cx} dx \\
&= \frac{(-1)^{k+l} k! l!}{c^{k+l+1-\alpha}} \int_0^{\infty} L_k^{\alpha}(t) L_l^{\alpha}(t) t^{\alpha} e^{-t} dt.
\end{align*}
Now, from the recursive relation for the monic generalized Laguerre polynomials $\lbrace \tilde{L}_k^{\alpha} \rbrace_{k \geq 0}$
\begin{equation}
\tilde{L}_{k+1}^{\alpha}(cx) = (cx-A_k) \tilde{L}_{k}^{\alpha}(cx)-B_k\tilde{L}_{k-1}^{\alpha}(cx),
\end{equation}
by (\ref{Lac}) we obtain
\begin{equation*}
c^{k+1} L_{k+1}^{\alpha,c}(x) = c \left( x-\frac{A_k}{c} \right) c^{k} L_{k}^{\alpha,c}(x)-B_k c^{k-1} L_{k-1}^{\alpha,c}(x),
\end{equation*}
and then (\ref{Lac_ric}).
\end{proof}

Using the polynomials $\lbrace L_k^{\alpha,c} \rbrace_{k \geq 0}$, the modified moments (see \ref{modified}) can be written as
\begin{equation*}
m_k = \int_0^{\infty} L_k^{\alpha,c} (x) x^{\alpha} e^{-cx} J_{\nu} (x)  dx + \int_0^{\infty} L_k^{\alpha,c} (x) x^{\alpha} e^{-cx} dx.
\end{equation*}
Clearly, by orthogonality, the second integral is zero for $k \geq 1$. Hence, for $k \geq 1$, by (\ref{Lac}) and the following explicit expression for the generalized Laguerre polynomials (see \cite[p.775]{Abramowitz})
\begin{equation*}
L_k^{\alpha}(x)= \sum_{j=0}^k (-1)^j \binom{k+\alpha}{k-j} \frac{1}{j!} x^j,
\end{equation*}
we have that
\begin{align}
m_k &= \int_0^{\infty} L_k^{\alpha,c} (x) x^{\alpha} e^{-cx} J_{\nu}(x) dx \notag \\
&= \frac{(-1)^k k!}{c^k} \int_0^{\infty} L_k^{\alpha} (cx) x^{\alpha} e^{-cx} J_{\nu}(x) dx \notag \\
&= \frac{(-1)^k k!}{c^k} \sum_{j=0}^k (-1)^j \binom{k+\alpha}{k-j} \frac{1}{j!} c^j \int_0^{\infty}  x^{\alpha+j} e^{-cx} J_{\nu}(x) dx \notag \\ 
&= \frac{(-1)^k k!}{c^k} \sum_{j=0}^k (-1)^j \binom{k+\alpha}{k-j} \frac{1}{j!} c^j \mu_{j,0}, \label{ricorsione_mod_mom}
\end{align}
where the last equality comes from (\ref{defCORE}).
Finally, for $k=0$ we obtain
\begin{equation*}
m_0 = \mu_{0,0} + \frac{\Gamma(\alpha+1)}{c^{\alpha+1}},
\end{equation*}
by \cite[sect.3.381, n.4]{GR}.

Using the modified moments, we can employ the Chebyshev algorithm, summarized in Algorithm \ref{ModChebyAlg}. We remark that the case $a_k = b_k = 0$ yields $p_k(x) = x^k$, and Algorithm \ref{ModChebyAlg} reduces to Algorithm \ref{ChebyAlg}.

\begin{algorithm} \label{ModChebyAlg}
Initialization
\begin{align*}
\alpha_0&=a_0 + \frac{m_1}{m_0}, \; \beta_0=m_0, \\
\sigma_{-1,l}&=0, \quad l=1,2,\ldots,2n-2, \\
\sigma_{0,l}&= m_l, \quad 0,1,\ldots,2n-1,
\end{align*}
for $k=1,2,\ldots,n-1$
\medskip

for $l=k,k+1,\ldots,2n-k-1$
\begin{align*}
\sigma_{k,l}&=\sigma_{k-1,l+1}-(\alpha_{k-1}-a_l) \sigma_{k-1,l}-\beta_{k-1} \sigma_{k-2,l}+b_l \sigma_{k-1,l-1}, \\
\alpha_k &= a_k + \frac{\sigma_{k,k+1}}{\sigma_{k,k}}-\frac{\sigma_{k-1,k}}{\sigma_{k-1,l-1}}, \quad \beta_k=\frac{\sigma_{k,k}}{\sigma_{k-1,k-1}}.
\end{align*}
\end{algorithm}

\begin{figure}
\begin{center}
\includegraphics[scale=0.35]{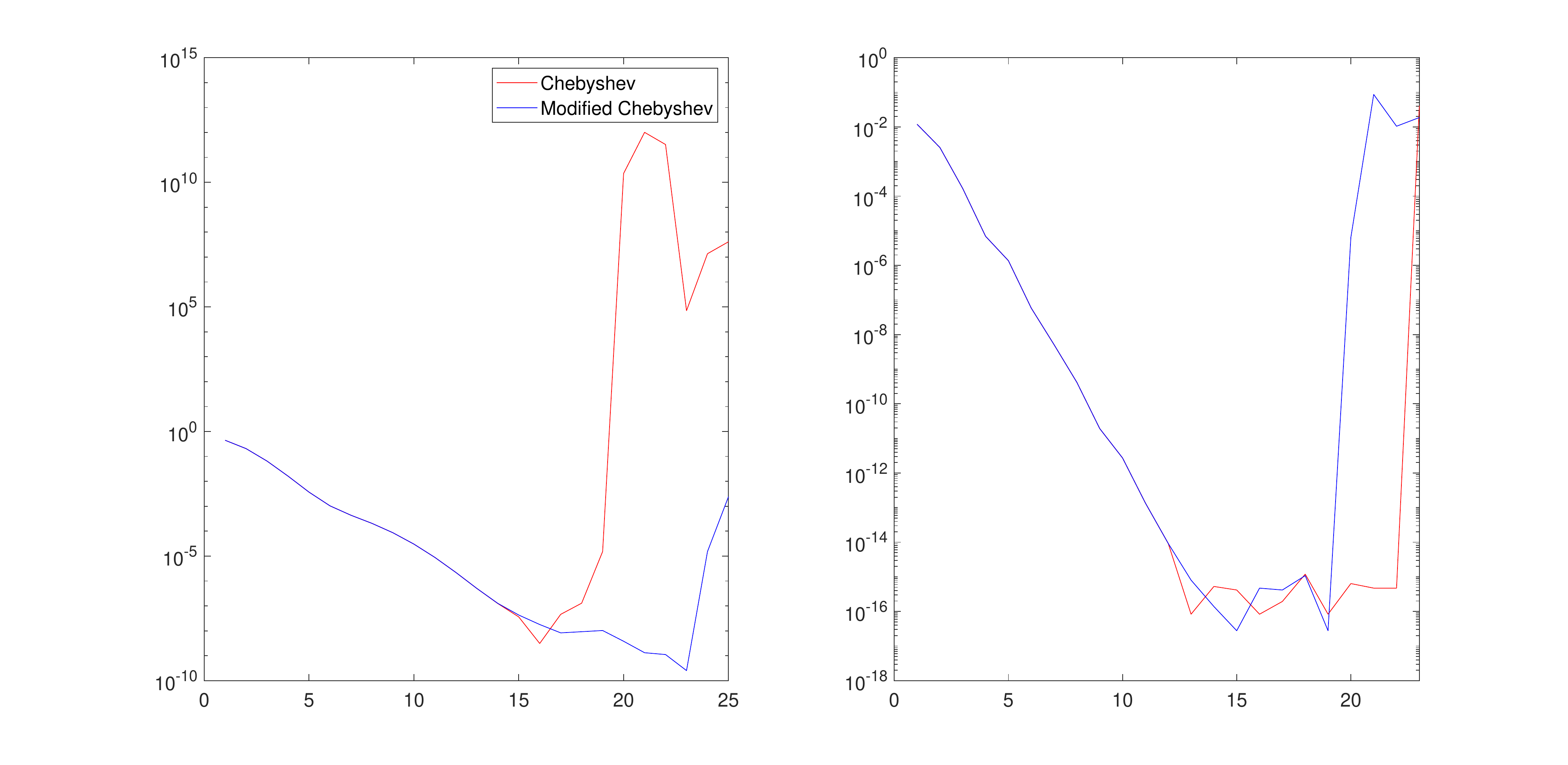}
\end{center}
\caption{Error behavior with respect to $n$ for $\nu =0.5$, $\alpha=0.5$, $c=0.2$ on the left and for $\nu =1$, $\alpha=0.5$, $c=0.7$ on the right. In both cases $f(x)=e^{-0.5x}$.}
\label{Figura2}
\end{figure}

In Figure \ref{Figura2} we compare the results of Algorithm \ref{ChebyAlg} and \ref{ModChebyAlg}. We provide only two representative examples that, nevertheless, are sufficient to say that Algorithm \ref{ModChebyAlg} in general allows to gain stability for further $5 \div 10$ iterations but in many cases there is no effective improvement.

\section{A preconditioned Cramer based approach} \label{section5}

Let
\[ M_k=
\begin{bmatrix}
\mu_0 & \mu_1 & \cdots & \mu_{k-1}\\
\mu_1 & \mu_2 & \cdots & \mu_k  \\
\vdots & \vdots &  & \vdots \\
\mu_{k-1} & \mu_k & \cdots &\mu_{2k-2}
\end{bmatrix} \in \mathbb{R}^{k \times k},
\]
be the moment matrix, and
\[ N_k=
\begin{bmatrix}
\mu_0 & \mu_1 & \cdots & \mu_{k-2} & \mu_k \\
\mu_1 & \mu_2 & \cdots & \mu_{k-1} & \mu_{k+1}  \\
\vdots & \vdots &  & \vdots & \vdots \\
\mu_{k-1} & \mu_k & \cdots &\mu_{2k-3} & \mu_{2k-1}
\end{bmatrix} \in \mathbb{R}^{k \times k}.
\]
It is known (see \cite{Davis}, \cite{Szego}) that the recurrence coefficients in (\ref{doppio_asterisco}) can also be written as 
\begin{equation} \label{ab_det}
\alpha_k = \frac{F_{k+1}}{D_{k+1}}-\frac{F_k}{D_k}, \quad \beta_k = \frac{D_{k-1}D_{k+1}}{D_k^2} \quad k \geq 0,
\end{equation}
where
\begin{align*}
D_k &=\det(M_k), \quad {\rm for} \quad k \geq 1, \\
F_k &= \det(N_k), \quad {\rm for} \quad k \geq 2,
\end{align*}
and
\begin{align*}
D_0&=D_{-1}=1,  \\
F_0&=0, \; F_1=\mu_1.
\end{align*}

Consider the linear system 
\begin{equation} \label{sistema1}
M_{k+1} x^{(k+1)} = e_{k+1},
\end{equation}
where $e_{k+1}=(0, \ldots, 0, 1)^T \in \mathbb{R}^{k+1}$.
In the following, we denote by $x_i^{(k+1)}$ the $i$-th component of the solution of (\ref{sistema1}).
First of all, we observe that, by Cramer's rule, 
\begin{equation*}
\frac{D_k}{D_{k+1}} = x_{k+1}^{(k+1)}.
\end{equation*}
Moreover, since
\begin{equation*}
\det(N_k) = -\det M_{k+1,(k)},
\end{equation*}
in which $ M_{k+1,(k)}$ is the matrix $ M_{k+1}$ with the $k$-th column substituted by the vector $e_{k+1}$, we have that
\begin{equation*}
\frac{F_k}{D_{k+1}} = -\frac{\det  M_{k+1,(k)}}{\det  M_{k+1}} = -x_k^{(k+1)}.
\end{equation*}
Hence, we obtain
\begin{equation*}
\frac{F_k}{D_k} = \frac{F_k}{D_{k+1}} \frac{D_{k+1}}{D_k}=-\frac{x_k^{(k+1)}}{x_{k+1}^{(k+1)}}.
\end{equation*}

In this setting, the coefficients $\alpha_k$ and $\beta_k$ can be expressed in terms of the components of the solutions of appropriate linear systems as follows:
\begin{equation} \label{alpha_beta_def}
\alpha_k = -\frac{x_{k+1}^{(k+2)}}{x_{k+2}^{(k+2)}} + \frac{x_{k}^{(k+1)}}{x_{k+1}^{(k+1)}}, \quad \beta_k = \frac{x_{k}^{(k)}}{x_{k+1}^{(k+1)}}, \quad k \geq 1,
\end{equation}
with
\begin{equation*}
\alpha_0 = \frac{\mu_1}{\mu_0}, \quad \beta_0 = \mu_0.
\end{equation*}

The system (\ref{sistema1}) rapidly becomes severely ill conditioned, so that the procedure does not offer any improvement with respect to the Chebyshev algorithm. 
Nevertheless, since $M_{k+1}$ is a symmetric positive definite matrix, the idea is to use a bilateral preconditioner in order to solve efficiently (\ref{sistema1}).
Analogously to the choice made for the modified Chebyshev approach, here we want to use as preconditioner the moment matrix corresponding to the generalized Laguerre polynomials.

Let $\eta_k$, $k \geq 0$, be the moments relative to the weight function $x^{\alpha}e^{-cx}$, given by
\begin{equation*}
\eta_k=\eta_k^{\alpha,c} = \int_0^{\infty} x^{k+\alpha} e^{-cx} dx = \frac{\Gamma(k+\alpha+1)}{c^{k+\alpha+1}},
\end{equation*}
where we have used again \cite[sect.3.381, n.4]{GR}. We can write
\begin{equation*}
\eta_k = \frac{1}{c^{\alpha+1}} \frac{\gamma_k}{c^k},
\end{equation*}
where
\begin{equation} \label{pallino}
\gamma_k=\gamma_k^{\alpha} = \int_0^{\infty} x^{k+\alpha} e^{-x} dx = \Gamma(k+\alpha+1)
\end{equation}
are the moments relative to the generalized Gauss-Laguerre rule. 
Hence, we can write the corresponding moment matrix
\[ M_k^{\alpha,c}= 
\begin{bmatrix}
\eta_0 & \eta_1 & \cdots & \eta_{k-1}\\
\vdots & \vdots &  & \vdots \\
\eta_{k-1} & \eta_k & \cdots & \eta_{2k-2}
\end{bmatrix} \in \mathbb{R}^{k \times k},
\]
as 
\begin{equation*}
M_k^{\alpha,c}=\frac{1}{c^{\alpha +1}} E_k M_k^{\alpha} E_k,
\end{equation*}
where 
\[ M_k^{\alpha}= 
\begin{bmatrix}
\gamma_0 & \gamma_1 & \cdots & \gamma_{k-1}\\
\vdots & \vdots &  & \vdots \\
\gamma_{k-1} & \gamma_k & \cdots & \gamma_{2k-2}
\end{bmatrix} \in \mathbb{R}^{k \times k},
\]
and $E_k = {\rm diag} \left( c^0,c^{-1}, \ldots, c^{1-k} \right)$. 

At this point, if we consider the Cholesky decomposition of $M_k^{\alpha}$
\begin{equation*}
M_k^{\alpha} = (R_k^{\alpha})^T R_k^{\alpha},
\end{equation*}
we have
\begin{equation*}
M_k^{\alpha,c} = \frac{1}{c^{\alpha+1}} E_k (R_k^{\alpha})^T R_k^{\alpha} E_k = (R_k^{\alpha,c})^T R_k^{\alpha,c},
\end{equation*}
where
\begin{equation*}
R_k^{\alpha,c} = \frac{1}{(\sqrt{c})^{\alpha+1}} R_k^{\alpha} E_k.
\end{equation*}
The following proposition provides the explicit expression for $R_k^{\alpha}$, and therefore for $R_k^{\alpha,c}$.

\begin{proposition}
The Cholesky decomposition of the matrix $M_k^{\alpha}$ is
\begin{equation*}
M_k^{\alpha} = (R_k^{\alpha})^T R_k^{\alpha},
\end{equation*}
with
\begin{equation*}
R^{\alpha}_{ij} = \frac{(j-1)!}{(j-i)!} \frac{\Gamma(\alpha+j)}{\sqrt{\Gamma(i)\Gamma(\alpha+i)}}, \quad {\rm for} \quad i \leq j \leq k.
\end{equation*}
\end{proposition}
\begin{proof}
Since the matrix $M_k^{\alpha}$ is symmetric, we can restrict the analysis to the case $i \leq j$. By (\ref{pallino}) we know that
\begin{equation*}
(M_k^{\alpha})_{ij}= \Gamma(i+\alpha+j-1).
\end{equation*}
Now,
\begin{align*}
((R_k^{\alpha})^T R_k^{\alpha})_{ij} &= \sum_{l=1}^{i} R_{lj}^{\alpha} R_{li}^{\alpha}  \\
&= \Gamma(j) \Gamma(\alpha+j)\Gamma(i) \Gamma(\alpha+i) \sum_{l=1}^{i} \frac{1}{(j-l)! (i-l)! \Gamma(l) \Gamma(\alpha+l)}.  
\end{align*}
Writing
\begin{equation} \label{relazione_bin}
\binom{x}{y} = \frac{\Gamma(x+1)}{\Gamma(y+1) \Gamma(x-y+1)},
\end{equation}
with $x= \alpha+i-1$ and $y=l-i$, we have
\begin{equation*}
\frac{1}{\Gamma(\alpha+l)(i-l)!} = \binom{\alpha+i-1}{i-l} \frac{1}{\Gamma(\alpha+i)}.
\end{equation*}
Moreover,
\begin{equation*}
\frac{1}{(j-l)!(l-1)!} = \binom{j-1}{l-1} \frac{1}{(j-1)!}.
\end{equation*}
Using the above relations, we obtain
\begin{equation*}
((R_k^{\alpha})^T R_k^{\alpha})_{ij}=\Gamma(\alpha+j)\Gamma(i) \sum_{l=1}^{i} \binom{\alpha+i-1}{i-l} \binom{j-1}{l-1}.
\end{equation*}
At this point, by the following slight modification of the Chu-Vandermonde identity (see \cite{Askey}),
\begin{equation*}
\sum_{u=1}^q \binom{t}{u-1} \binom{s-t}{q-u} = \binom{s}{q-1},
\end{equation*}
we obtain
\begin{equation*}
((R_k^{\alpha})^T R_k^{\alpha})_{ij}= \Gamma(\alpha+j)\Gamma(i) \binom{\alpha+j+i-2}{i-1}.
\end{equation*}
Using again (\ref{relazione_bin}), with $x=\alpha+j+i-2$ and $y=i+1$, it holds
\begin{equation*}
\binom{\alpha+j+i-2}{i-1}= \frac{\Gamma(i+j+\alpha-1)}{\Gamma(\alpha+j)\Gamma(i)},
\end{equation*}
and finally
\begin{equation*}
((R_k^{\alpha})^T R_k^{\alpha})_{ij} = \Gamma(i+\alpha+j-1).
\end{equation*}
\end{proof}

We observe that the matrix $R_k^{\alpha}$ can be written as
\begin{equation*}
R_k^{\alpha}=D_k \tilde{R}_k^{\alpha},
\end{equation*}
with
\begin{equation*}
(\tilde{R}_k^{\alpha})_{ij} = \frac{(j-1)! \Gamma(\alpha+j)}{(j-i)! \Gamma(i) \Gamma(\alpha+i)}, \quad {\rm for} \quad i \leq j,
\end{equation*}
and $D_k$ diagonal matrix such that 
\begin{equation*}
(D_k)_{ii}=\sqrt{\Gamma(i) \Gamma(\alpha+i)}.
\end{equation*}
Since
\begin{equation*}
(R_k^{\alpha})^{-1} = (\tilde{R}_k^{\alpha})^{-1} (D_k^{\alpha})^{-1}.
\end{equation*}
and
\begin{equation*}
(\tilde{R}_k^{\alpha})_{ij}^{-1} = (-1)^{i+j} (\tilde{R}_k^{\alpha})_{ij},
\end{equation*}
we have that the explicit expression for $(R_k^{\alpha})^{-1}$ is given by
\begin{equation*}
(R_k^{\alpha})_{ij}^{-1} = (-1)^{i+j} \frac{\sqrt{(j-1)! \Gamma(\alpha+j)}}{(j-i)! \Gamma(i) \Gamma(\alpha+i)}, \quad {\rm for} \quad i \leq j. 
\end{equation*}
Therefore, the matrix $(R_k^{\alpha,c})^{-1}$ can be written as
\begin{equation*}
(R_k^{\alpha,c})^{-1} = (\sqrt{c})^{\alpha+1} E_k^{-1} (\tilde{R}_k^{\alpha})^{-1} (D_k^{\alpha})^{-1},
\end{equation*}
with
\begin{equation} \label{matrice IRac}
(R_k^{\alpha,c})^{-1}_{ij}=\frac{\sqrt{c}^{\alpha+1}c^{i-1}(-1)^{i+j}\sqrt{(j-1)! \Gamma(\alpha+j)}}{(j-i)! \Gamma(i) \Gamma(\alpha+i) }.
\end{equation}

Finally, the linear system (\ref{sistema1}) can be preconditioned as
\begin{equation} \label{sistema2}
(R_{k+1}^{\alpha,c})^{-T} M_{k+1} (R_{k+1}^{\alpha,c})^{-1} y^{(k+1)} = (R_{k+1}^{\alpha,c})^{-T} e_{k+1},
\end{equation}
with
\begin{equation} \label{soluzione}
x^{(k+1)}=(R_{k+1}^{\alpha,c})^{-1} y^{(k+1)}.
\end{equation}
Since the matrix $M_{k+1}$ can be written as $M_{k+1}=M_{k+1}^{\alpha,c}+M_{k+1,0}^{\alpha,c}$, where $M_{k+1,0}^{\alpha,c}$ is the matrix of the core moments defined by equation (\ref{defCORE}), we have that
\begin{align}
(R_{k+1}^{\alpha,c})^{-T} M_{k+1} (R_{k+1}^{\alpha,c})^{-1} &= (R_{k+1}^{\alpha,c})^{-T}(M_{k+1}^{\alpha,c}+M_{k+1,0}^{\alpha,c}) (R_{k+1}^{\alpha,c})^{-1} \notag \\
&= (R_{k+1}^{\alpha,c})^{-T} M_{k+1}^{\alpha,c} (R_{k+1}^{\alpha,c})^{-1}+(R_{k+1}^{\alpha,c})^{-T}M_{k+1,0}^{\alpha,c}(R_{k+1}^{\alpha,c})^{-1} \notag \\
&=  I_{k+1}+(R_{k+1}^{\alpha,c})^{-T}M_{k+1,0}^{\alpha,c}(R_{k+1}^{\alpha,c})^{-1} := Q_{k+1}, \label{Q_k}
\end{align}
where $I_{k+1}$ is the identity matrix. The system (\ref{sistema2}) becomes
\begin{equation*}
(I_{k+1}+ (R_{k+1}^{\alpha,c})^{-T} M_{k+1,0}^{\alpha,c} (R_{k+1}^{\alpha,c})^{-1}) y^{(k+1)} = (R_{k+1}^{\alpha,c})^{-T} e_{k+1}.
\end{equation*}
In Table \ref{tabella_condizionamento} we show the remarkable effect of the preconditioning. 
\begin{table}[h]
\[
\begin{array}{cccccccc}
\toprule
k  & 5  &  10 & 15 & 20 &  25 &  30  \\ 
\midrule
\kappa_2 (M_k) & 2.4e+13 & 7.8e+32 & 1.0e+51 & 2.2e+69 &  4.6e+88 & 1.1e+107 \\
\kappa_2 (Q_k) & 1.0e+00 & 1.3e+00 & 1.4e+00 & 1.4e+00 & 1.5e+00 & 1.6e+00 \\
\bottomrule &  &  &   
\end{array}
\]
\caption{The Euclidean condition number of the matrix $M_k$ and of the preconditioned matrix $Q_k$, defined in (\ref{Q_k}), for different values of $k$. In this example $\nu=0.9$, $\alpha=0.1$ and $c=0.1$.}
\label{tabella_condizionamento}
\end{table}

We observe that, since $(R_{k+1}^{\alpha,c})^{-1}$ is an upper triangular matrix, the components of the solution used in (\ref{ab_det}) can be written as
\begin{align}
x_k^{(k+1)} &= (R_{k+1}^{\alpha,c})^{-1}_{kk} y_k^{(k+1)}+(R_{k+1}^{\alpha,c})^{-1}_{k,k+1} y_{k+1}^{(k+1)}, \label{remark1} \\
x_{k+1}^{(k+1)} &= (R_{k+1}^{\alpha,c})^{-1}_{k+1,k+1} y_{k+1}^{(k+1)}. \label{remark2}
\end{align}

We notice that the numerical implementation of the procedure to calculate $\alpha_k$ and $\beta_k$ as in (\ref{alpha_beta_def}), by using expressions (\ref{remark1}) and (\ref{remark2}), starts to show instability around $k=60 \div 70$, depending on the parameters, when $x_k^{(k+1)}$ and $x_{k+1}^{(k+1)}$ are close to the underflow.
In order to gain more stability the idea is to rewrite the coefficients $\alpha_k$ and $\beta_k$,  for $k \geq 1$, in terms of the components of the vectors $y^{(k)},y^{(k+1)},y^{(k+2)}$, defined in (\ref{soluzione}), and exploit the relation (\ref{matrice IRac}). Indeed, we observe that for $i \sim j \sim k$,
\begin{equation*}
(R_k^{\alpha,c})_{ij}^{-1} \sim \frac{c^k}{k!},
\end{equation*}
and therefore $y_i^{(k+1)} \gg x_i^{(k+1)}$ for $i=k,k+1$. By (\ref{ab_det}) and (\ref{matrice IRac}), we obtain
\begin{align}
\alpha_k &= -\frac{\sqrt{(k+1)(\alpha+k+1)}}{c} \left(\frac{y_{k+1}^{(k+2)}}{y_{k+2}^{(k+2)}} -\sqrt{(k+1)(\alpha+k+1)} \right) \notag \\
&+ \frac{\sqrt{k(\alpha+k)}}{c} \left(\frac{y_{k}^{(k+1)}}{y_{k+1}^{(k+1)}} -\sqrt{k(\alpha+k)} \right), \label{alphak30} \\
\beta_k &= \frac{\sqrt{k(\alpha+k)}}{c} \left( \frac{y_k^{(k)}}{y_{k+1}^{(k+1)}} \right). \label{betak31}
\end{align}

The final procedure, explained in Algorithm \ref{algoritmo3}, allows to work with $80 \div 90$ points, dependently on the parameters.

\begin{algorithm} \label{algoritmo3}
Define $\alpha_0$, $\beta_0$, $ \beta_1$, $y^{(1)}$, $y^{(2)}$.

\smallskip
for $k=2, \ldots, n-1$

\smallskip
\quad calculate $y^{(k+1)}$ by solving (\ref{sistema2})

\smallskip
\quad $\beta_k \leftarrow y^{(k+1)}$, $y^{(k)}$ by (\ref{betak31})

\smallskip
\quad $\alpha_{k-1} \leftarrow y^{(k+1)}$, $y^{(k)}$ by (\ref{alphak30})

\smallskip
end
\end{algorithm}

\begin{figure}
\begin{center}
\includegraphics[scale=0.35]{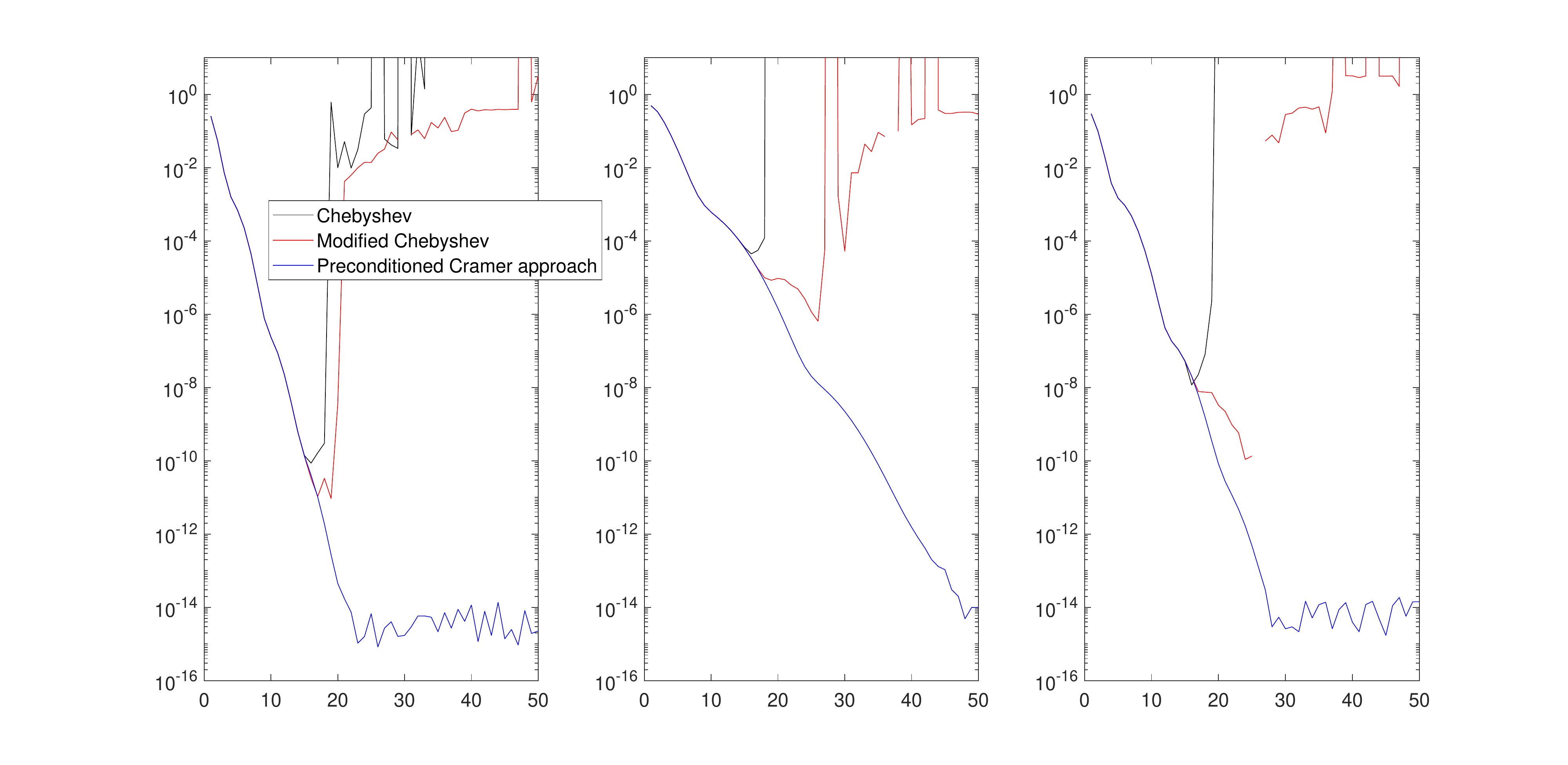}
\caption{Error histories for $\nu =1$, $\alpha=0.7$, $c=0.3$ on the left, for $\nu =0.9$, $\alpha=0.1$, $c=0.1$ in the middle and for $\nu =1.5$, $\alpha=0.5$, $c=0.2$ on the right. In all cases $f(x)=e^{-0.5x}$.}  \label{figura3}
\end{center}
\end{figure}

\begin{figure} 
\begin{center}
\includegraphics[scale=0.35]{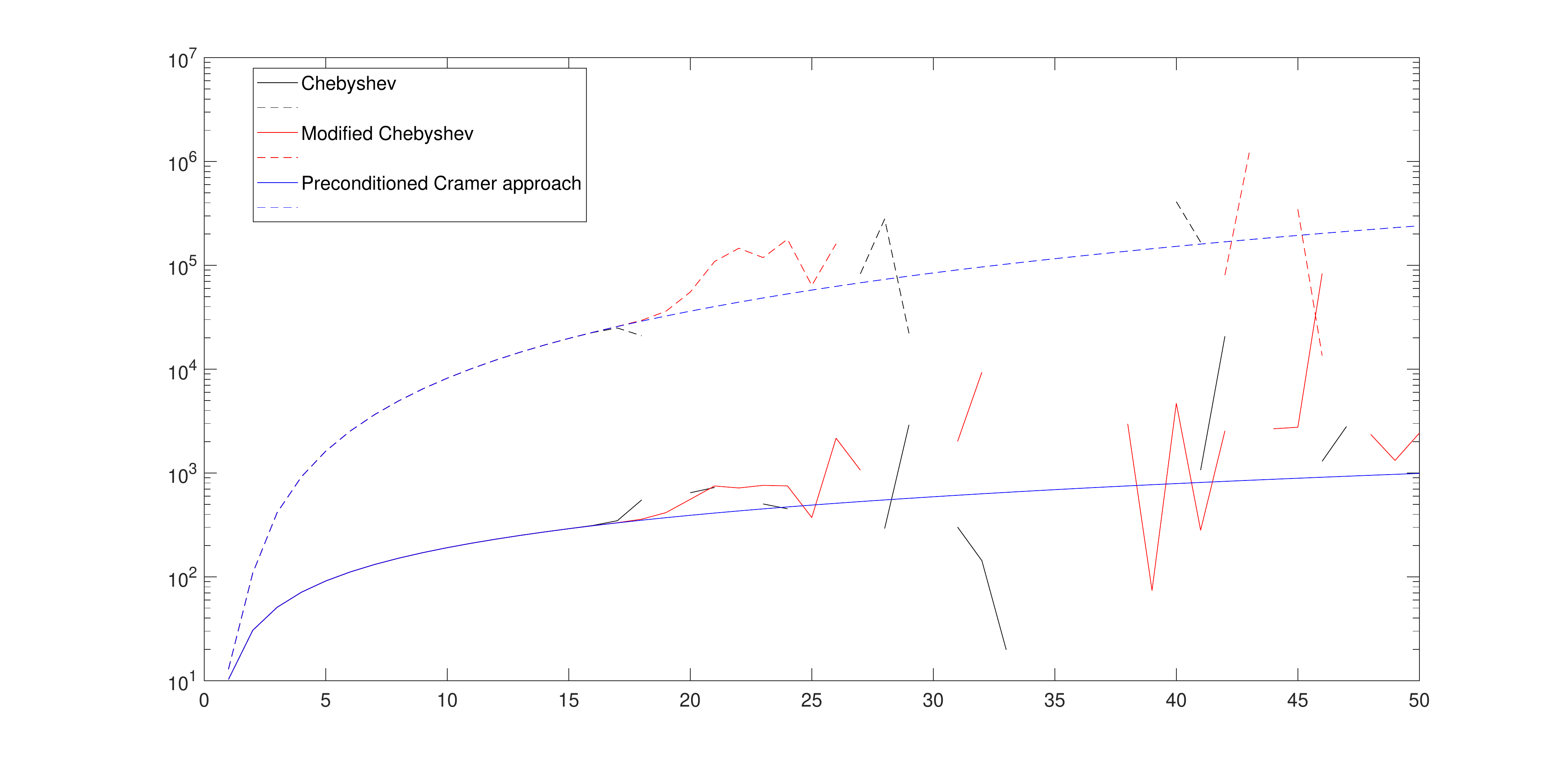}
\caption{Plot of the coefficients $\alpha_k$ (solid lines) and $\beta_k$ (dashed lines) for $\nu=0.9$, $\alpha=0.1$ and $c=0.1$.} \label{figura4}
\end{center}
\end{figure}

In Figure \ref{figura3} we compare the results of Algorithm \ref{ChebyAlg}, \ref{ModChebyAlg} and \ref{algoritmo3}. For all the examples we can say that only the preconditioned Cramer based approach allows to achieve an absolute error around the machine precision, while the Chebyshev and modified Chebyshev algorithms lose stability much earlier. In fact, as shown in Figure \ref{figura4}, Algorithm \ref{ChebyAlg} and \ref{ModChebyAlg} start to provide inaccurate values of the coefficients $\alpha_k$ and $\beta_k$ for $k$ around $15 \div 25$, while Algorithm \ref{algoritmo3} is definitely more stable. Since the plot is in logarithmic scale, the missing parts of the curves are relative to negative entries.

\section{Electromagnetic fields} \label{section6}

In this section we deal with an interesting application arising in geophysical electromagnetic (EM) survey.
We consider the theoretical EM response, i.e. the electromagnetic fields components, over a $N$-layered earth due to vertical magnetic dipoles above the surface, composed of a transmitter coil and couples of receiver coils. The receiver couples are placed at different distances (offsets) from the transmitter coil. In this case, the electromagnetic induction effect, encoded in the first-order linear differential equations, produces eddy alterning currents in the soil which on their turn, induce response EM fields. Under the assumption that each layer is characterized by a certain conductivity $\sigma_i$ and thickness $h_i$ (the deeper layer is assumed to have infinite thickness), the general integral solutions of Maxwell equations are given by (see \cite{CR1} and \cite{WH})
\begin{align*}
H_{z}^{(N)}& =\frac{m}{4\pi }\int_{0}^{\infty }(1+R_{0}(\lambda )e^{-2H \lambda})\lambda^{2}J_{0}(\lambda r)d\lambda , \\
H_{\rho }^{(N)}& =\frac{m}{4\pi }\int_{0}^{\infty }(1-R_{0}(\lambda
)e^{-2H \lambda})\lambda ^{2}J_{1}(\lambda r)d\lambda , 
\end{align*}
where $m$ is the magnetic moment, $H$ is the height of the dipole with respect to the surface and $r$ is the offset. 
In the above formulas $
R_{0}(\lambda )$ is the reflection term, recursively defined by 
\begin{equation*}
\begin{split}
R_{0}(\lambda )& =\frac{R_{1}(\lambda )+\Psi _{1}(\lambda )}{R_{1}(\lambda
)\Psi _{1}(\lambda )+1}, \\
R_{j}(\lambda )& =\frac{R_{j+1}(\lambda )+\Psi _{j+1}(\lambda )}{%
R_{j+1}(\lambda )\Psi _{j+1}(\lambda )+1}e^{-2u_{j}(\lambda )h_{j}},\quad
j=1,...,N-1, \\
R_{N}(\lambda )& =0,
\end{split}
\end{equation*}
with
\begin{equation*}
\Psi _{j}(\lambda )=\frac{u_{j-1}(\lambda )-u_{j}(\lambda )}{u_{j-1}(\lambda
)+u_{j}(\lambda )},\quad j=1,...,N,
\end{equation*}%
in which $u_{0}(\lambda )=\lambda $ and $u_{j}(\lambda )=\sqrt{\lambda^{2}-k_{j}^{2}}$, $k_{j}=\sqrt{-i\omega \mu \sigma _{j}}$, for $j=1,\ldots ,N$, where $\omega$ is the angular frequency and $\mu$ is the magnetic permeability of vacuum.
We refer to \cite[Section 4]{WH} and the reference therein for an exhaustive background.
Since in the case of conductivity of geological materials only the imaginary part of the fields are considered, using the change of variable $\lambda r=x$, we obtain
\begin{align}
\Im (H_z^{(N)}) &= \frac{m}{4\pi r^3}\int_{0}^{\infty }\Im \left( R_{0}\left( \frac{x}{r} \right) \right) e^{-\frac{2H}{r} x} x^{2}J_{0}(x)dx, \label{Hz}\\
\Im (H_{\rho}^{(N)}) &= -\frac{m}{4\pi r^3}\int_{0}^{\infty }\Im \left( R_{0}\left( \frac{x}{r} \right) \right) e^{-\frac{2H}{r} x} x^{2}J_{1}(x)dx. \label{Hrho}
\end{align}

In the numerical experiments we use Algorithm \ref{algoritmo3} to evaluate the fields (\ref{Hz}) and (\ref{Hrho}) in the case of a $3$-layered underground model. Referring to (\ref{integrale}), in our examples we set $\nu=0,1$, $\alpha=0$, $c=\frac{2H}{r}$ and $f(x)=\Im \left( R_{0}\left( \frac{x}{r} \right) \right) x^2$.
Regarding the choice of the parameters $\sigma_i$ and $h_i$, i.e. of the underground models, we consider real life values of river levees (see e.g. \cite{argini}).

\begin{figure}
\begin{center}
\includegraphics[scale=0.35]{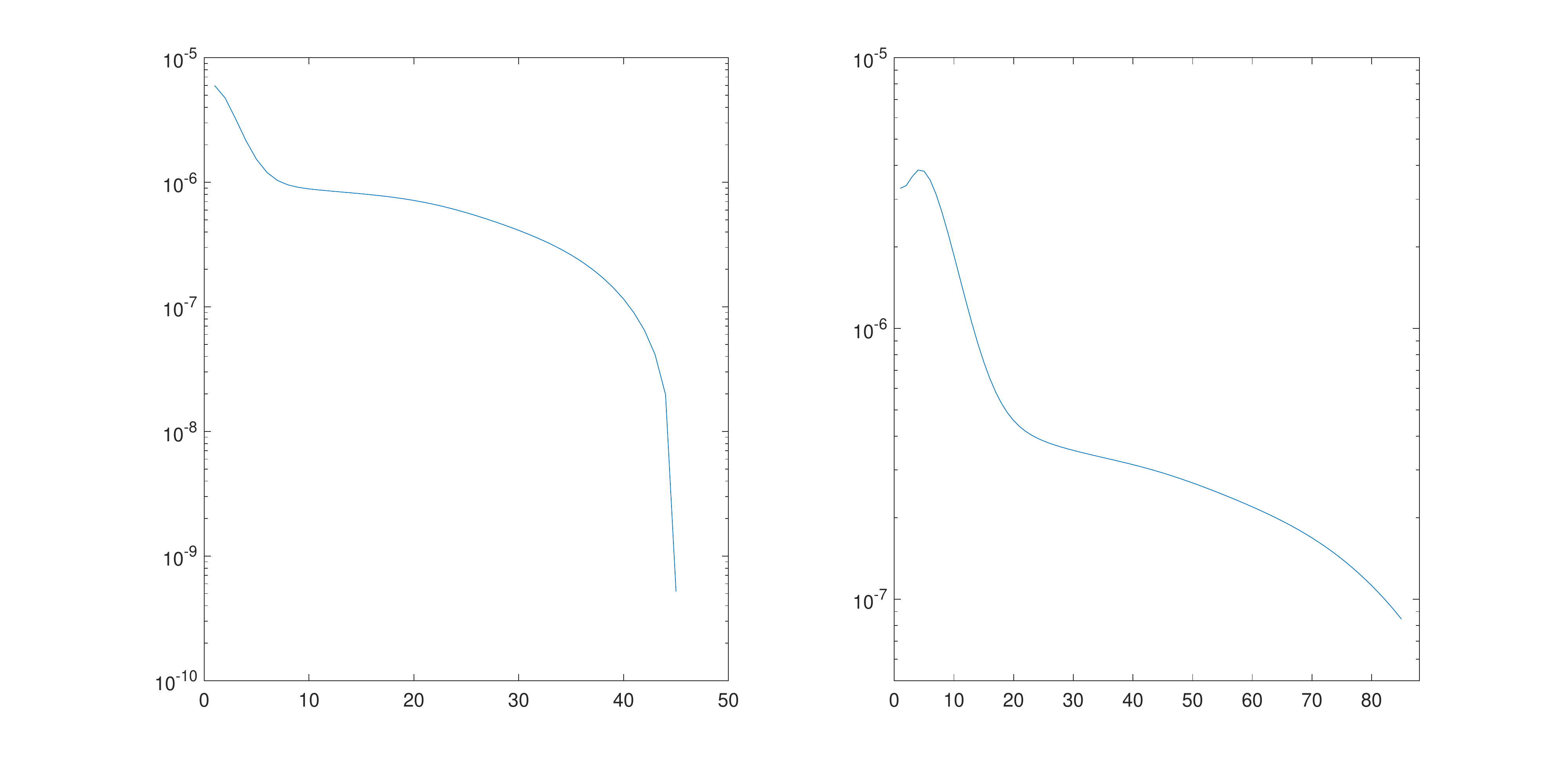}
\end{center}
\caption{Error history for the computation of $\Im (H_z^{(3)})$ with respect to $n$ for parameters $H=0.4$ $m$, $r=8$ $m$, $h_1=2.5$ $m$, $h_1=0.5$ $m$, $\sigma_1=0.05$ $S/m$, $\sigma_2=0.0049$ $S/m$, $\sigma_3=0.0182$ $S/m$ on the left and for parameters $H=0.2$ $m$, $r=8$ $m$, $h_1=2.5$ $m$, $h_1=0.5$ $m$, $\sigma_1=0.033$ $S/m$, $\sigma_2=0.1$ $S/m$, $\sigma_3=0.01$ $S/m$ on the right.} \label{figuraHz}
\end{figure}

\begin{figure}
\begin{center}
\includegraphics[scale=0.35]{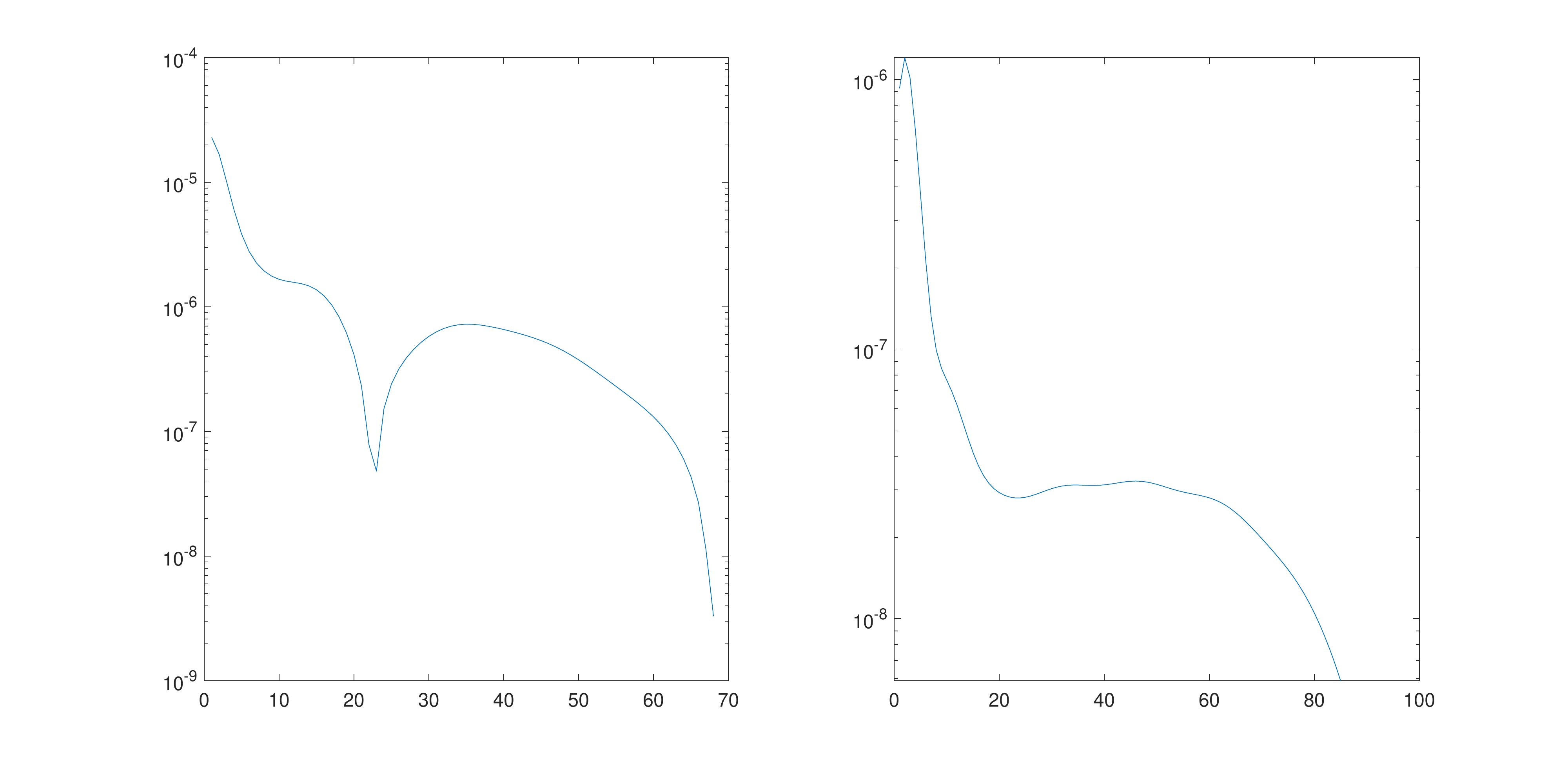}
\end{center}
\caption{Error history for the computation of $\Im (H_{\rho}^{(3)})$ with respect to $n$ for parameters $H=0.4$ $m$, $r=8$ $m$, $h_1=2.5$ $m$, $h_1=0.5$ $m$, $\sigma_1=0.333$ $S/m$, $\sigma_2=0.02$ $S/m$, $\sigma_3=0.1$ $S/m$ on the left and for parameters $H=0.4$ $m$, $r=8$ $m$, $h_1=2.5$ $m$, $h_1=0.5$ $m$, $\sigma_1=0.033$ $S/m$, $\sigma_2=0.1$ $S/m$, $\sigma_3=0.01$ $S/m$ on the right.} \label{figuraHrho}
\end{figure}

In Figure \ref{figuraHz} and \ref{figuraHrho} we provide the absolute error between the approximated fields $\Im (H_z^{(3)})$ and $\Im (H_{\rho}^{(3)})$, and a corresponding reference solution (see e.g. \cite{G-S}, \cite{Emdpler}). In all examples we stop the procedure when the error is less than $10^{-8}$. The main reason of this choice is that for these parameters the method works with $k_{{\rm max}} \sim 85$, but in order to reach the machine precision more points are necessary. To overcome this issue, the symbolic computation and quadruple-precision arithmetic can be adopted (see e.g. \cite{Gautschi_simbolico}), but this is beyond the purpose of the present paper. 

\section{Conclusions}

In this work a Gaussian type quadrature rule for the computation of integrals involving fractional powers, exponentials and Bessel functions of the first kind, is presented. In this framework, the techniques commonly used in the computation of the coefficients of the three-term recurrence relation, for the corresponding orthogonal polynomials, are the standard and the modified Chebyshev algorithm. Since it is well known that the results of these methods can be inaccurate for growing number of quadrature points and especially for unbounded intervals of integrations, an alternative and very stable approach, based on the preconditioning of the moment matrix, is developed. In particular, an algorithm, which exploits the Cramer rule to compute the coefficients by solving a linear system with the moment matrix, is presented. The numerical experiments confirm the reliability of this preconditioned Cramer based approach and shows that it is definitely more stable than the modified Chebyshev algorithm, since it allows to work with further $40 \div 60$ points, depending on the parameters.

We remark that, in principle, the approach can be applied to each weight function that is not so far to the standard ones, because it is necessary to be able to construct the preconditioner. Finally, we also point out that, similarly to the Gauss Laguerre rule, the weights decay exponentially and therefore a truncated approach can also be introduced as well.

\section*{Acknowledgements}

This work was partially supported by GNCS-INdAM, FRA-University of Trieste and CINECA under HPC-TRES program award number 2019-04. The authors are members of the INdAM research group GNCS.

\end{document}